\documentclass[reqno]{amsart}
\usepackage{amsmath, amsthm, amssymb, amstext}

\usepackage[left=3cm,right=3cm,top=3cm,bottom=3cm]{geometry}
\usepackage{hyperref,xcolor}
\hypersetup{
 pdfborder={0 0 0},
 colorlinks,
}
\usepackage{enumitem}
\newtheorem{theorem}{Theorem}
\newtheorem{remark}[theorem]{Remark}

\newtheorem{proposition}[theorem]{Proposition}
\newtheorem{corollary}[theorem]{Corollary}

\newtheorem{example}[theorem]{Example}

\DeclareMathOperator*{\divergenz}{div}              %
\DeclareMathOperator*{\esssup}{ess ~sup}         %

\newcommand{\N}{\mathbb{N}}

\newcommand{\R}{\mathbb{R}}

\newcommand{\RN}{\mathbb{R}^N}

\newcommand{\Om}{\Omega}
\newcommand{\rand}{\partial\Omega}

\newcommand{\into}{\int_{\Omega}}

\renewcommand{\l}{\left}
\renewcommand{\r}{\right}
\numberwithin{theorem}{section}
\numberwithin{equation}{section}

\title[Existence and $L^{\infty}$-estimates for elliptic equations involving
convolution]
{Existence and $L^{\infty}$-estimates for elliptic equations involving
convolution}

\author[G.\,Marino]{Greta Marino}
\address[G.\,Marino]{Technische Universit$\ddot{\text a}$t Chemnitz, Fakult\"at f\"ur Mathematik,
Reichenhainer Stra\ss e 41, 09126 Chemnitz, Germany}
\email{greta.marino@mathematik.tu-chemnitz.de}

\author[D.\,Motreanu]{Dumitru Motreanu}
\address[D.\,Motreanu]{D\'epartement de Math\'emathiques, Universit\'e de Perpignan, 66860 Perpignan, France}
\email{motreanu@univ-perp.fr}

\subjclass[2010]{35J60, 35B45, 35J25, 44A35}

\keywords{Moser iteration, boundedness of solutions, elliptic
operators of divergence type, critical growth on the boundary,
convolution}

\begin{document}

\begin{abstract}
In this paper, with a fixed $p\in (1,+\infty)$ and a bounded domain
$\Omega \subset \mathbb{R}^N$, $N \ge 2$, whose boundary $\partial\Omega$
fulfills the Lipschitz regularity, we study the following boundary
value problem
    \begin{equation}
    \label{prob} \tag{P}
    \begin{aligned}
    - \divergenz \mathcal A(x, u, \nabla u)+a|u|^{p-2}u&= \mathcal B(x, \rho
\ast E(u),\nabla(\rho\ast E(u))) \quad && \text{in } \Omega, \\
    \mathcal A(x, u, \nabla u) \cdot \nu&= \mathcal C(x, u) && \text{on
} \partial\Omega,
    \end{aligned}
    \end{equation}
where $\mathcal A\colon \Omega\times\mathbb{R} \times \mathbb{R}^N
\to \mathbb{R}^N$, $\mathcal B\colon \Omega \times \mathbb{R} \times
\mathbb{R}^N \to \mathbb{R}$, $\mathcal C\colon \partial\Omega
\times \mathbb{R} \to \mathbb{R}$ are Carath\'eodory functions,
$a>0$ is a constant, $E\colon W^{1,p}(\Omega)\to
W^{1,p}(\mathbb{R}^N)$ is an extension operator related to $\Omega$,
and $\rho$ is an integrable function on $\mathbb{R}^N$. This is a
novel problem that involves the nonlocal operator assigning to $u$
the convolution $\rho \ast E(u)$ of $\rho$ with $E(u)$. Under
verifiable conditions, we prove the existence of a (weak) solution
to problem \eqref{prob} by using the surjectivity theorem for
pseudomonotone operators. Moreover, through a modified version of
Moser iteration up to the boundary initiated in \cite{MW, MW1} we
show that (any) weak solution to \eqref{prob} is bounded.
\end{abstract}

\maketitle

\section{Introduction}\label{S1}

Let $\Omega \subset \mathbb{R}^N$, $N \ge 2$, be a bounded domain with a
Lipschitz continuous boundary $\partial\Omega$ and let $p\in
(1,+\infty)$ be a real number. It is well-known that there exists an
extension operator $E\colon W^{1,p}(\Omega)\to
W^{1,p}(\mathbb{R}^N)$ meaning that $E$ is a linear map satisfying
$$
E(u)|_\Omega=u,\quad \forall \, u\in W^{1,p}(\Omega)
$$
and for which there exists a constant $C= C(\Omega)>0$ depending
only on $\Omega$ such that
    \[
    \begin{split}
    \|E(u)\|_{L^p(\mathbb{R}^N)}&\leq C(\Omega)\|u\|_{L^p(\Omega)} \\
    \text{and} \quad  \|E(u)\|_{W^{1,p}(\mathbb{R^N})}&\leq C(\Omega) \|u\|_{W^{1,p}(\Omega)},
    \end{split}
     \quad \forall \, u\in W^{1,p}(\Omega)
    \]
(see \cite {A,B}). In the terminology of \cite{A} such a map $E$ is
called a $(1,p)$-extension operator for $\Omega$. Generally, the
extension operators are constructed by using reflection maps and
partitions of unity. For the rest of the paper, we fix an extension
operator $E\colon W^{1,p}(\Omega)\to W^{1,p}(\mathbb{R}^N)$.

We state the following boundary value problem
\begin{equation}\label{problem}
\begin{aligned}
- \divergenz \mathcal A(x, u, \nabla u)+a|u|^{p-2}u&= \mathcal B(x,
\rho
\ast E(u),\nabla(\rho\ast E(u))) \quad && \text{in } \Omega, \\
\mathcal A(x, u, \nabla u) \cdot \nu&= \mathcal C(x, u) && \text{on
} \partial\Omega,
\end{aligned}
\end{equation}
where $a>0$ is a constant, $\nu(x) $ denotes the outer unit normal
of $\Omega$ at $x\in \rand$, $\rho \ast E(u)$ stands for the
convolution product of some integrable function $\rho$ on
$\mathbb{R}^N$ with $E(u)$, and $\mathcal A$, $\mathcal B$, $
\mathcal C $ are Carath\'eodory functions satisfying suitable
$p$-structure growth conditions. Due to the presence of convolution,
problem \eqref{problem} is nonlocal. Furthermore, in the statement
of problem \eqref{problem} we have full dependence on the solution
$u$ and on its gradient $\nabla u$, which makes the problem highly
non-variational, so the variational methods are not applicable. The
boundary condition in \eqref{problem} is nonhomogeneous and includes
the Robin boundary condition.

The starting point of this work has been the elliptic problem in
\cite{MM} with homogeneous  Dirichlet boundary condition
\begin{equation}\label{problem1}
\begin{aligned}
-\Delta_{p}u-\mu\Delta_{q}u&=f(x,\rho\ast u,\nabla(\rho\ast u))\ && \mbox{in $\Omega$}, \\
u&=0 && \mbox{on } \partial\Omega,
\end{aligned}
\end{equation}
involving the $p$-Laplacian $\Delta_{p}$ and the $q$-Laplacian
$\Delta_{q}$ with $1<q<p<+\infty$, where for the first time the
boundary value problem with convolution for solution and its
gradient was considered. Any solution $u\in W^{1,p}_0(\Omega)$ of
\eqref{problem1} can be identified with $E(u)\in
W^{1,p}(\mathbb{R^N})$ obtained by extension with zero outside
$\Omega$. In this case both $\rho$ and $u$ are integrable functions
on $\mathbb{R}^N$ and the convolution $\rho\ast u$ in
\eqref{problem1} makes sense. This is no longer possible for
\eqref{problem} because we have $u\in W^{1,p}(\Omega)$ and the
extension by zero outside $\Omega$ generally does nor produce an
element of $W^{1,p}(\mathbb{R}^N)$. Here is the essential point
where the extension operator $E $ is necessary in \eqref{problem}.

Finally, among papers involving quasilinear elliptic equations with convection term we can refer to \cite{Mar-Wink}.

The aim of this paper is two-fold: to establish an existence result
for \eqref{problem} and to provide a priori estimates for the
solutions to \eqref{problem} up to the boundary showing their
uniform boundedness. The proof of existence of solutions to
\eqref{problem} relies on the theory of pseudomonotone operators and
properties of convolution and extension operator. In order to prove
a priori estimates for problem \eqref{problem} and show the
boundedness of its solutions we develop a modified version of Moser
iteration originating in \cite{MW, MW1}.

First of all we recall that the critical exponents corresponding to $p$ in
$\Omega$ and on $\partial\Omega$ are denoted by $p^*$ and $p_*$,
respectively (see Section \ref{S2}). 

For the existence result, our assumptions are as follows.

\begin{itemize}

\item[(A)] The maps $\mathcal A\colon \Omega \times \mathbb{R} \times \mathbb{R}^N \to \mathbb{R}^N$,
$\mathcal B\colon \Omega \times\mathbb{R}\times\mathbb{R}^N \to
\mathbb{R}$, and $\mathcal C\colon \partial\Omega
\times\mathbb{R}\to \mathbb{R}$ are Carath\'eodory functions (i.e.,
they are measurable in the first variable and continuous in the
others) satisfying the following conditions:
\[
\begin{aligned}
& \text{(A1)} \quad \vert \mathcal A(x, s, \xi) \vert \le a_1
\vert\xi \vert^{p-1}
+ a_2 \vert s \vert^{p-1}+ a_3 \qquad && \text{for a.e. } x \in \Omega, \\
& \text{(A2)} \quad \mathcal A(x, s, \xi-\xi') \cdot (\xi-\xi')>0&& \text{for a.e. } x \in \Omega, \\
& \text{(A3)} \quad \mathcal A(x, s, \xi) \cdot \xi \ge a_4 \vert
\xi \vert^p- a_5 && \text{for a.e. } x \in \Omega, \\
& \text{(A4)} \quad \vert \mathcal B(x, s, \xi) \vert \le f(x)
+ b_1 \vert s \vert^{\alpha_1}+ b_2 \vert \xi \vert^{\alpha_2} && \text{for a.e. } x \in \Omega, \\
& \text{(A5)} \quad \vert \mathcal C(x, s) \vert \le c_1 \vert s
\vert^{\alpha_3}+ c_2 && \text{for a.e. } x \in \rand,
\end{aligned}
\]
for all $s \in \mathbb{R}$ and $\xi,\xi'\in\mathbb{R}^N$,
$\xi\not=\xi'$, with positive constants $a_i, b_j, c_k \, (i \in
\{1, \dots, 5\}, j, k \in \{1,2\}) $, with
\begin{eqnarray}\label{21}
\alpha_1,\alpha_2,\alpha_3\in[0,p-1)
\end{eqnarray}
and a nonnegative function $f\in L^{r'}(\Omega)$ with $r\in [1,
p^*)$. 
\end{itemize}

Assumptions (A1)-(A2) are the Leray-Lions conditions, while (A3) is
a coercivity condition. In problem \eqref{problem1} we have
$\mathcal A(x,s,\xi)=|\xi|^{p-2}\xi+\mu|\xi|^{q-2}\xi$, with
$1<q<p<+\infty$ and $\mu\geq 0$, which fulfills these assumptions.
The maps $\mathcal B$ and $\mathcal C$ are only subject to the
growth conditions (A4)-(A5).

By a (weak) solution to problem \eqref{problem} we mean any function
$u\in W^{1,p}(\Omega)$ verifying
\begin{equation}\label{weak-sol}
\int_{\Omega} \mathcal A(x, u, \nabla u) \cdot \nabla \varphi
dx+a\int_{\Omega}|u|^{p-2}u\varphi dx= \int_{\Omega} \mathcal B(x,
\rho \ast E(u), \nabla(\rho \ast E(u))) \varphi dx+
\int_{\partial\Omega} \mathcal C(x,u)\varphi d\sigma
\end{equation}
for all $\varphi\in W^{1,p}(\Omega)$. Under assumptions (A), all the
integrals in \eqref{weak-sol} are finite for $u,\varphi \in
W^{1,p}(\Omega)$, thus the definition of weak solution is
meaningful. In the same spirit, $u \in W^{1,p}_0(\Omega)$ is a
(weak) solution to \eqref{problem1} if
\[
\int_{\Omega} \l(\vert \nabla u \vert^{p-2} + \mu \vert \nabla u
\vert^{q-2} \r) \nabla u \cdot \nabla \varphi dx = \int_{\Omega}
\mathcal B(x, \rho \ast u, \nabla (\rho \ast u)) \varphi dx
\]
holds for every $ \varphi \in W^{1,p}_0(\Omega)$.

\begin{theorem}\label{T1}
Let $\Omega\subset\mathbb{R}^N$ be a bounded domain with a Lipschitz
continuous boundary $\partial\Omega$ endowed with the extension
operator $E\colon W^{1,p}(\Omega)\to W^{1,p}(\mathbb{R}^N)$ and let
$\rho\in L^1(\mathbb{R}^N)$. If hypotheses \emph{(A)} are satisfied,
then there exists a (weak) solution to problem \eqref{problem}.
\end{theorem}

The proof of Theorem \ref{T1} is the object of Section \ref{S3}.

Now we turn to the uniform boundedness of solutions to problem
\eqref{problem}. We formulate the following
hypotheses.
\begin{itemize}

\item[(H)] The maps $\mathcal A\colon \Omega \times \R \times \RN \to \RN$,
$\mathcal B\colon \Omega \times \R \times \RN \to \R $ and $
\mathcal C\colon \rand \times \R \to \R$ are Carath\'eodory
functions satisfying the conditions
\[
\begin{aligned}
& \text{(H1)} \quad \vert \mathcal A(x, s, \xi) \vert \le a_1
\vert\xi \vert^{p-1}
+ a_2 \vert s \vert^{p^* \frac{p-1}{p}}+ a_3 \qquad && \text{for a.e. } x \in \Omega, \\
& \text{(H2)} \quad \mathcal A(x, s, \xi) \cdot \xi \ge a_4 \vert
\xi \vert^p- a_5
\vert s \vert^{p^*}- a_6 && \text{for a.e. } x \in \Omega, \\
& \text{(H3)} \quad \vert \mathcal B(x, s, \xi) \vert \le  f(x)
+ b_1 \vert s \vert^{\alpha_1}+ b_2 \vert \xi \vert^{\alpha_2} && \text{for a.e. } x \in \Omega, \\
& \text{(H4)} \quad \vert \mathcal C(x, s) \vert \le c_1 \vert s
\vert^{p_*-1}+ c_2 && \text{for a.e. } x \in \rand,
\end{aligned}
\]
for all $s \in \mathbb{R}$ and $\xi \in \mathbb{R}^N$, with
nonnegative constants $a_i, b_j, c_k \, (i \in \{1, \dots, 6\}, j, k
\in \{1,2\}) $ and $\alpha_1, \alpha_2 $ such that
\begin{equation}\label{alpha}
0 \le \alpha_1< p^*-p, \qquad 0 \le \alpha_2< \min \left\{p-1,
\frac{p}{p^*}(p^*-p)\right\},
\end{equation}
and a nonnegative function $ f \in L^{r'}(\Omega)$, with $ r \in [1,
p^*/p)$.
\end{itemize}

\begin{theorem}\label{T2}
Let $\Omega\subset \RN $ be a bounded domain with a Lipschitz
continuous boundary $\partial\Omega$ endowed with the extension
operator $E\colon W^{1,p}(\Omega)\to W^{1,p}(\mathbb{R}^N)$ and let
$\rho\in L^1(\mathbb{R}^N)$. Assume that hypotheses \emph{(H)} are
satisfied. Then, every (weak) solution $u\in W^{1,p}(\Omega) $ to
problem \eqref{problem} belongs to $L^\infty(\Omega)$ with the trace
$\gamma u\in L^{\infty}(\partial\Omega)$.
\end{theorem}

The proof of Theorem \ref{T2} is  given in Section \ref{S4}.

Combining Theorems \ref{T1} and \ref{T2} we obtain the following
existence result of bounded solutions to problem \eqref{problem}.

\begin{corollary}\label{C1}
Let $\Omega \subset \RN$ be a bounded domain with a Lipschitz
continuous boundary $\partial\Omega$ endowed with the extension
operator $E\colon W^{1,p}(\Omega)\to W^{1,p}(\mathbb{R}^N)$ and let
$\rho \in L^1(\mathbb{R}^N)$. Assume that hypotheses
\emph{(A1)-(A3)}, \emph{(A4)} with $\alpha_2$ as in \eqref{alpha},
and \emph{(A5)} are satisfied. Then, there exists a (weak) solution
$u\in W^{1,p}(\Omega) $ to problem \eqref{problem} which belongs to
$L^\infty(\Omega)$ and whose trace $\gamma u$ is an element of
$L^{\infty}(\partial\Omega)$.
\end{corollary}

Corollary \ref{C1} is a direct consequence of Theorems \ref{T1} and
\ref{T2} noticing that Theorems \ref{T1} and \ref{T2} can be
simultaneously applied.

We illustrate the applicability of our results by an example using
the extension operator constructed in \cite[page 275]{B}.

\begin{example}\label{E1}
Consider in $\mathbb{R}^2$ the rectangular domains
$\Omega=(0,1)\times (0,1)$, $\Omega_1=(0,1)\times (-1,1)$,
$\Omega_2=(-1,1)\times (-1,1)$, $\Omega_3=(-1,1)\times (-1,3)$,
$\tilde{\Omega}=(-1,3)\times (-1,3)$. We introduce the maps
$R_1\colon W^{1,p}(\Omega)\to W^{1,p}(\Omega_1)$, $R_2\colon
W^{1,p}(\Omega_1)\to W^{1,p}(\Omega_2)$, $R_3\colon
W^{1,p}(\Omega_2)\to W^{1,p}(\Omega_3)$, and $R_4\colon
W^{1,p}(\Omega_3)\to W^{1,p}(\tilde{\Omega})$, respectively, by
\begin{align*}
(R_1u)(x_1,x_2)=
\begin{cases}
u(x_1,x_2) \qquad & \text{if } x_2>0, \\
u(x_1,-x_2) & \text{if } x_2<0
\end{cases}
\end{align*}
for all  $u\in W^{1,p}(\Omega)$ and $(x_1,x_2)\in \Omega$,
\begin{align*}
(R_2u)(x_1,x_2)=
\begin{cases}
u(x_1,x_2) \qquad & \text{if } x_1>0, \\
u(-x_1,x_2) & \text{if } x_1<0
\end{cases}
\end{align*}
for all  $u\in W^{1,p}(\Omega_1)$ and $(x_1,x_2)\in \Omega_1$,
\begin{align*}
(R_3u)(x_1,x_2)=
\begin{cases}
u(x_1,x_2) \qquad & \text{if } x_2<1, \\
u(x_1,2-x_2) & \text{if } x_2>1
\end{cases}
\end{align*}
for all  $u\in W^{1,p}(\Omega_2)$ and $(x_1,x_2)\in \Omega_2$, and
\begin{align*}
(R_4u)(x_1,x_2)=
\begin{cases}
u(x_1,x_2) \qquad & \text{if } x_1<1, \\
u(2-x_1,x_2) & \text{if } x_1>1
\end{cases}
\end{align*}
for all  $u\in W^{1,p}(\Omega_3)$ and $(x_1,x_2)\in \Omega_3$.

For a fixed $\psi\in C^1(\tilde{\Omega})$ with $\psi=1$ on $\Omega$
and ${\rm supp}\ \psi\subset \tilde{\Omega}$, the linear map
$E\colon W^{1,p}(\Omega)\to W^{1,p}(\mathbb{R}^2)$ which carries
each $u\in W^{1,p}(\Omega)$ to the function $Eu\in
W^{1,p}(\mathbb{R}^2)$ obtained by extending $\psi(R_4\circ
R_3\circ R_2\circ R_1 u)$ with zero outside $\tilde{\Omega}$ is an
extension operator. Accordingly, given a constant $a>0$, a function $\rho\in
L^1(\mathbb{R}^2)$ and a Carath\'eodory function $\mathcal B\colon
\Omega \times \mathbb{R}\times \mathbb{R}^2\to \mathbb{R}$
satisfying \emph{(H)} and \eqref{alpha} we state the Neumann problem
\begin{equation*}
\begin{aligned}
- \Delta_pu+a\vert u\vert^{p-2}u&=\mathcal B(x, \rho
\ast E(u),\nabla(\rho\ast E(u))) \quad && \text{in } \Omega, \\
\vert \nabla u \vert^{p-2}\nabla u\cdot \nu&= 0 && \text{on }
\partial\Omega.
\end{aligned}
\end{equation*}
A frequent form of $\mathcal B$ is $\mathcal
B(x,s,\xi)=g(s)+h(\xi)$. Our  results apply to the stated problem.
\end{example}

The rest of the paper is organized as follows. Section \ref{S2}
contains preliminaries to be used in the sequel. In Section \ref{S3}
we prove Theorem \ref{T1}. In Section \ref{S4} we prove Theorem
\ref{T2}.

\section{Preliminaries}\label{S2}

The Euclidean norm of $ \RN$ is denoted by $ \vert \cdot \vert $,
while the notation $ \cdot $ stands for the standard inner product
on $ \RN$. By $ \vert \cdot \vert $ we also denote the Lebesgue
measure on $\RN$. In the rest of the paper, for every $r\in (1,
+\infty)$ we denote by $r'$ its H\"older conjugate, that is $ r'=
\frac{r}{r-1}$.

For any $r\in [1, \infty)$ and a domain $\Omega\subset
\mathbb{R}^N$, we denote by $L^r(\Omega)$ and $W^{1,r}(\Omega)$ the
usual Lebesgue and Sobolev spaces equipped with the norms
\begin{align}\label{norms}
\begin{aligned}
\Vert u \Vert_{L^r(\Omega)}&= \l(\int_\Omega \vert u \vert^r dx \r)^{\frac{1}{r}},  \\
\Vert u \Vert_{W^{1,r}(\Omega)}&= \l(\int_\Omega \vert \nabla u
\vert^r dx \r)^{\frac{1}{r}}+ \l(\int_\Omega \vert u \vert^r dx
\r)^{\frac{1}{r}}.
\end{aligned}
\end{align}
Recall that the norm of $L^{\infty}(\Omega)$ is
\begin{align*}
\Vert u \Vert_{L^{\infty}(\Omega)}= \esssup_{\Omega} \vert u \vert.
\end{align*}
For any $u \in W^{1,r}(\Omega)$, set $u^{\pm}:= \max\{\pm u, 0\}$,
which yields
\begin{align}
\label{pm} u^{\pm} \in W^{1,r}(\Omega), \quad \vert u \vert= u^+ +
u^-, \quad u= u^+- u^-.
\end{align}
By the Sobolev embedding theorem there exists a linear continuous
embedding $i\colon W^{1,r}(\Omega) \to L^{r^*}(\Omega)$, where the
corresponding critical exponent $r^*$ in the domain is given by
\begin{align*}
\label{critical-dom}
r^*=
\begin{cases}
\frac{Nr}{N-r} \qquad & \text{if } r< N, \\
+\infty & \text{if } r \ge N.
\end{cases}
\end{align*}
The boundary $\partial\Omega$ is endowed with the
$(N-1)$-dimensional Hausdorff (surface) measure. The measure of
$\partial\Omega$ is denoted by $|\partial\Omega|$. The Lebesgue spaces
$ L^s(\partial\Omega)$, with $1\le s \le +\infty$, have the norms
\begin{align*}
\Vert u \Vert_{L^{s}(\partial\Omega)}= \l(\int_{\partial\Omega}
\vert u \vert^s d\sigma \r)^{\frac{1}{s}} \quad (1 \le s<+\infty),
\qquad \Vert u \Vert_{L^{\infty}(\partial\Omega)}=
\esssup_{\partial\Omega} \vert u \vert.
\end{align*}
There exists a unique linear continuous map $\gamma\colon
W^{1,r}(\Omega) \to L^{r_*}(\partial\Omega)$, called the trace map,
characterized by $\gamma(u)= u|_{\partial\Omega}$ whenever $u\in
W^{1,r}(\Omega) \cap C(\overline{\Omega})$, where $r_*$ is the
corresponding critical exponent on the boundary defined as
\begin{align*}
r_*=
\begin{cases}
\frac{(N-1)r}{N-r} \qquad & \text{if } r< N, \\
+\infty & \text{if } r \ge N.
\end{cases}
\end{align*}
As usual, the subspace of $W^{1,r}(\Omega)$ consisting of zero trace
elements is denoted $W^{1,r}_0(\Omega)$. For the sake of notational
simplicity, we drop the use of the symbol $ \gamma$ writing simply
$u$ in place of $\gamma u$. We refer to \cite{A} for the theory of
Sobolev spaces.

The following propositions are useful in the proof of our
boundedness result.

\begin{proposition}(\cite[Proposition 2.2]{MW})\label{prop1}
Let  $u \in L^p(\Omega)$, $1< p< +\infty$, be
nonnegative. If it holds
\[
\Vert u \Vert_{L^{\alpha_n}(\Omega)} \le C
\]
for a constant $C> 0$ and a sequence $(\alpha_n) \subset
\mathbb{R}_+$ such that $ \alpha_n \to +\infty $ as $n \to \infty$,
then $ u \in L^{\infty}(\Omega)$.

\end{proposition}

\begin{proposition}(\cite[Proposition 2.4]{MW})
\label{prop3} 
Let $1< p< +\infty$ and let  $u \in
W^{1,p}(\Omega) \cap L^{\infty}(\Omega)$. Then,  $u\in
L^{\infty}(\rand)$.

\end{proposition}

Recall that for $ \rho\in L^1(\mathbb{R}^N)$ and $u\in
W^{1,p}(\mathbb{R}^N)$, with $1<p<+\infty$, the convolution $\rho
\ast u$ is defined by
\[
(\rho \ast u)(x):= \int_{\RN} \rho(x- y) u(y) dy \quad \text{for
a.e. } x \in \mathbb{R}^N.
\]
The weak partial derivatives of the convolution $\rho \ast u$ are
expressed by
\begin{equation*}
\frac{\partial}{\partial x_i}(\rho \ast u)= \rho \ast \frac{\partial
u} {\partial x_i} \quad \text{for}\ i= 1, \dots, N.
\end{equation*}
Thanks to Tonelli's and Fubini's theorems as well as H\"older's
inequality, there hold
\begin{equation*}
\Vert \rho \ast u \Vert_{L^r(\mathbb{R}^N)}\le \Vert \rho
\Vert_{L^1(\mathbb{R}^N)} \Vert u \Vert_{L^r(\mathbb{R}^N)}
\end{equation*}
for every $r \in [1, p^*]$ and
\begin{equation}
\label{holder2} \l\Vert \rho \ast \frac{\partial u}{\partial x_i}
\r\Vert_{L^p(\mathbb{R}^N)} \le \Vert \rho \Vert_{L^1(\mathbb{R}^N)}
\l\Vert\frac{\partial u}{\partial x_i} \r\Vert_{L^p(\mathbb{R}^N)}
\quad \text{for}\ i= 1, \dots, N
\end{equation}
(see  \cite[Theorem 4.15]{B}). Taking into account the fact that the
function $t \mapsto t^{1/2}$ is sublinear as well as the function
$t\mapsto t^p$ is convex on $(0,+\infty)$ and \eqref{holder2}, it
follows that
    \begin{equation*}
    \begin{split}
    \Vert \nabla(\rho\ast u)\Vert_{L^p(\mathbb{R}^N)}^p &=
\int_{\mathbb{R}^N} \vert \nabla(\rho \ast u) \vert^p dx=
\int_{\mathbb{R}^N} \l(\sum_{i=1}^N \l
(\rho \ast \frac{\partial u}{\partial x_i} \r)^2 \r)^{p/2} dx \\
    & \le \int_{\mathbb{R}^N} \l(\sum_{i=1}^N \l\vert \rho \ast
\frac{\partial u}{\partial x_i} \r\vert \r)^p dx \le N^{p-1}
\int_{\mathbb{R}^N} \sum_{i=1}^N \l
\vert \rho \ast \frac{\partial u}{\partial x_i} \r\vert^p dx \\
    &\leq N^p \Vert \rho \Vert_{L^1(\mathbb{R}^N)}^p \Vert \nabla u
\Vert_{L^p(\mathbb{R}^N)}^p.
\end{split}
\end{equation*}
Finally, we recall the main theorem on the pseudomonotone operators
that will be used to prove our existence result. Let $X$ be a
reflexive Banach space endowed with the norm $\|\cdot\|$. The norm
convergence is denoted by $\rightarrow$ and the weak convergence by
$\rightharpoonup$. We denote by $X^\ast$ the topological dual of $X$
and by $\left\langle \cdot \,, \cdot\right\rangle$ the duality
pairing between $X$ and $X^\ast$. A map $A\colon X\rightarrow
X^\ast$ is called bounded if it maps bounded sets to bounded sets.
It is said to be coercive if there holds
\begin{eqnarray*}
\lim_{\|u\|\rightarrow +\infty}\frac{\left\langle Au,u\right\rangle
}{\|u\|}=+\infty.
\end{eqnarray*}
Finally, $A$ is called pseudomonotone if $u_n \rightharpoonup u$ in
$X$ and
\begin{eqnarray*}
\limsup_{n\rightarrow +\infty}\left\langle
Au_n,u_n-u\right\rangle\leq 0
\end{eqnarray*}
imply
\begin{eqnarray*}
\left\langle Au,u-w\right\rangle \leq \liminf_{n\rightarrow
+\infty}\left\langle Au_n,u_n-w\right\rangle, \quad \forall \, w\in
X.
\end{eqnarray*}
The surjectivity theorem for pseudomonotone operators reads as
follows (see, e.g., \cite{CLM}).
\begin{theorem}\label{TMPMO}
Let $X$ be a reflexive Banach space, let $A\colon X\rightarrow
X^\ast$ be a pseudomonotone, bounded and coercive operator, and let
$g\in X^\ast$. Then, there exists at least a solution $u\in X$ to
the equation $Au=g$.
\end{theorem}

\section{Proof of Theorem \ref{T1}}\label{S3}

Throughout the proof of the theorem, we will denote by $C_i$, $i \in
\N$, constants which depend  on the given data.

With a fixed $\rho\in L^1(\mathbb{R}^N)$ and an extension operator
$E\colon W^{1,p}(\Omega)\to W^{1,p}(\mathbb{R}^N)$, we introduce the
nonlinear operator $T\colon W^{1,p}(\Omega)\rightarrow
(W^{1,p}(\Omega))^*$ by
\begin{align}\label{operator T}
\langle Tu,\varphi\rangle&=\int_{\Omega} \mathcal A(x, u, \nabla u)
\cdot \nabla \varphi dx+a\int_{\Omega}|u|^{p-2}u\varphi dx\\
\nonumber & \quad - \int_{\Omega} \mathcal B(x,\rho\ast E(u),
\nabla(\rho \ast E(u))) \varphi dx- \int_{\partial\Omega} \mathcal
C(x,u)\varphi d\sigma
\end{align}
for all $u,\varphi\in W^{1,p}(\Omega)$. Assumption (A) guarantees
that  $T$ is well defined.

Let us show that $T$ is also bounded. Indeed, fix $\varphi \in
W^{1,p}(\Omega) $ such that  $ \Vert \varphi \Vert_{W^{1,p}(\Omega)}
\le 1$. Then,
\begin{equation}
\label{bound1}
\begin{split}
\vert\langle Tu,\varphi\rangle \vert &\leq \int_{\Omega}
\vert\mathcal A(x, u, \nabla u)\vert \vert
\nabla \varphi\vert dx+a\int_{\Omega}\vert u\vert^{p-1}\vert\varphi\vert dx \\
& \quad+ \int_{\Omega} \vert\mathcal B(x,\rho\ast E(u),\nabla(\rho
\ast E(u)))\vert \vert \varphi\vert
dx+\int_{\partial\Omega}\vert\mathcal C(x,u)\vert\vert\varphi\vert
d\sigma.
\end{split}
\end{equation}
We estimate the terms of the inequality above separately.  First of
all, observe that
    \begin{equation}
    \label{bound2}
    \begin{split}
    \into \vert \mathcal A(x, u, \nabla u) \vert \vert \nabla \varphi \vert dx & \le \into \l(a_1 \vert \nabla u \vert^{p-1}+ a_2 \vert u \vert^{p-1}+ a_3 \r) \vert \nabla \varphi \vert dx \\
    & \le a_1 \Vert \nabla u \Vert_{L^p(\Omega)}^{p-1} \Vert \nabla \varphi \Vert_{L^p(\Omega)}+ a_2 \Vert u \Vert_{L^p(\Omega)}^{p-1} \Vert \nabla \varphi \Vert_{L^p(\Omega)}\\
    & \quad + a_3 \vert \Omega \vert^{\frac{p-1}{p}} \Vert \nabla \varphi \Vert_{L^p(\Omega)} \\
    & \le a_1 \Vert \nabla u \Vert_{L^p(\Omega)}^{p-1}+ a_2 \Vert u \Vert_{L^p(\Omega)}^{p-1}+ C_1,
    \end{split}
    \end{equation}
as well as
    \begin{equation}
    \label{bound3}
    a \into \vert u \vert^{p-1} \vert \varphi \vert dx \le a \Vert u \Vert_{L^p(\Omega)}^{p-1} \Vert \varphi \Vert_{L^p(\Omega)} \le a \Vert u \Vert_{L^p(\Omega)}^{p-1}.
    \end{equation}
Thanks to (A4) we also have
    \begin{equation}
    \label{bound4}
    \begin{split}
    &\into \vert \mathcal B(x, \rho \ast E(u), \nabla(\rho \ast E(u))) \vert \vert \varphi \vert dx  \\
    & \le \into \l(f(x)+ b_1 \vert \rho \ast E(u) \vert^{\alpha_1}+ b_2 \vert \nabla(\rho \ast E(u)) \vert^{\alpha_2} \r) \vert \varphi \vert dx.
    \end{split}
    \end{equation}
We consider the terms in \eqref{bound4} separately. First note that
H\"older's inequality gives
    \begin{equation}
    \label{bound-f}
    \begin{split}
    \into f(x) \vert \varphi \vert dx &\le \Vert f \Vert_{L^{r'}(\Omega)} \Vert \varphi \Vert_{L^r(\Omega)}  \\
    & \le \Vert f \Vert_{L^{r'}(\Omega)} \Vert \varphi \Vert_{L^p(\Omega)} |\Omega|^{\frac{p-r}{pr}} \\
    & \le C_2.
    \end{split}
    \end{equation}
Moreover, exploiting the properties of $E$ and of the convolution
and the Sobolev embedding we have
    \begin{equation}
    \label{conv1}
    \begin{split}
    b_1 \into |\rho \ast E(u) |^{\alpha_1}  |\varphi | dx & \le \Vert \rho \ast E(u) \Vert_{L^{p^*}(\RN)}^{\alpha_1} \Vert \varphi \Vert_{L^{\frac{p^*}{p^*-\alpha_1}}(\Omega)} \\
    & \le C_3 \Vert \rho \Vert_{L^1(\RN)}^{\alpha_1} \Vert E(u) \Vert_{L^{p^*}(\RN)}^{\alpha_1} \Vert \varphi \Vert_{L^{p^*}(\Omega)} \\
    & \le C_4 \Vert \rho \Vert_{L^1(\RN)}^{\alpha_1} \Vert u \Vert_{L^{p^*}(\Omega)}^{\alpha_1} \Vert \varphi \Vert_{W^{1,p}(\Omega)} \\
    & \le C_5 \Vert u \Vert_{W^{1,p}(\Omega)}^{\alpha_1},
    \end{split}
    \end{equation}
as well as
    \begin{equation}
    \label{conv2}
    \begin{split}
    b_2 \into \vert \nabla (\rho \ast E(u))\vert^{\alpha_2} \vert \varphi \vert dx & \le b_2 \Vert \nabla(\rho \ast E(u))\Vert_{L^p(\RN)}^{\alpha_2} \Vert \varphi \Vert_{L^{\frac{p}{p- \alpha_2}}(\Omega)} \\
    & \le C_6 \Vert \rho \Vert_{L^1(\RN)}^{\alpha_2} \Vert \nabla E(u) \Vert_{L^p(\RN)}^{\alpha_2} \Vert \varphi \Vert_{W^{1,p}(\Omega)} \\
    & \le C_7 \Vert \nabla u \Vert_{L^p(\RN)}^{\alpha_2} \le C_7 \Vert u \Vert_{W^{1,p}(\Omega)}^{\alpha_2}.
    \end{split}
    \end{equation}
Finally, hypothesis (A5) gives the following estimate for the
boundary term in \eqref{bound1}
    \begin{equation}
    \label{bound5}
    \begin{split}
    \int_{\rand} |\mathcal C(x, u)| |\varphi| d\sigma & \le
    \int_{\rand} \l(c_1 \vert u \vert^{\alpha_3}+ c_2 \r) \vert \varphi \vert d\sigma \\
     & \le c_1 \Vert u \Vert_{L^{p_*}(\rand)}^{\alpha_3} \Vert \varphi \Vert_{L^{\frac{p_*}{p_*- \alpha_3}}(\rand)}+ c_2 \vert \rand \vert^{\frac{p-1}{p}} \Vert \varphi \Vert_{L^p(\rand)} \\
    & \le c_1 \Vert u \Vert_{W^{1,p}(\Omega)}^{\alpha_3}+ C_8.
    \end{split}
    \end{equation}
Taking into account \eqref{bound2}-\eqref{bound5} and applying once
again the Sobolev embedding, from \eqref{bound1} we derive
    \[
    \vert \langle Tu, \varphi \rangle \vert \le C_9(\|u\|_{W^{1,p}(\Omega)}^{\beta}+1),
    \]
for all $\Vert \varphi \Vert_{W^{1,p}(\Omega)} \le 1$, with $
\beta:= \max\{p-1, \alpha_1, \alpha_2, \alpha_3\}$. This in turn
implies
    \[
    \Vert Tu \Vert_{(W^{1,p}(\Omega))^*} \le C_9(\|u\|_{W^{1,p}(\Omega)}^{\beta}+1),
    \]
which shows that $T$ is bounded.

Now we prove that  $T$  is pseudomonotone. Toward this, let
$(u_n)_{n \in \N} \subset W^{1,p}(\Omega)$ be a sequence satisfying
$u_n\rightharpoonup u$ for some $u\in W^{1,p}(\Omega)$ and
    \begin{equation}\label{12}
    \limsup_{n\to +\infty} \, \langle Tu_n,u_n-u\rangle \leq 0.
    \end{equation}
By H\"older's inequality and Rellich-Kondrachov compact embedding
theorem it follows that, passing to a subsequence if necessary,
    \begin{equation}
    \label{15}
    \begin{split}
    \left|\int_{\Omega}|u_n|^{p-2}u_n(u_n-u) dx\right| & \leq \int_\Omega |u_n|^{p-1}|u_n-u| dx \\
    &\leq \|u_n\|_{L^p(\Omega)}^{p-1} \|u_n-u\|_{L^p(\Omega)} \to
0\ \quad \mbox{as}\ n\to\infty.
    \end{split}
    \end{equation}
With a similar argument already exploited in
\eqref{bound-f}-\eqref{conv2} we have
    \[
    \begin{split}
    &\left|\int_{\Omega} \mathcal B(x,\rho\ast E(u_n),\nabla(\rho \ast
E(u_n)))(u_n-u)dx\right|\\
    & \quad \leq \int_{\Omega} \l\vert\mathcal B(x,\rho\ast E(u_n),\nabla(\rho
\ast E(u_n))) \r\vert \vert
u_n-u\vert dx \\
    &\quad \leq C_{10} \|u_n-u\|_{L^r(\Omega)}+ C_{11} \|u_n\|_{W^{1,p}(\Omega)}^{\alpha_1} \|u_n-u\|_{L^{\frac{p^*}{p^*-\alpha_1}}(\Omega)} \\
    & \quad \quad + C_{12} \|u_n\|_{W^{1,p}(\Omega)}^{\alpha_2}
\|u_n-u\|_{L^{\frac{p}{p-\alpha_2}}(\Omega)}
    \end{split}
    \]
for all $u\in W^{1,p}(\Omega)$. Since
$$
r, \frac{p^*}{p^*-\alpha_1}, \frac{p}{p-\alpha_2}<p^*,
$$
we can apply  the Rellich-Kondrachov compact embedding theorem to
the previous estimate, which gives
\begin{eqnarray}\label{16}
\lim_{n\to+\infty} \int_{\Omega} \mathcal B(x,\rho\ast
E(u_n),\nabla(\rho \ast E(u_n)))(u_n-u)dx=0.
\end{eqnarray}
Finally,  assumption (A), H\"older's inequality and the compactness
of the trace mappings due to the inequalities
$$
p, \frac{p_*}{p_*-\alpha_3}<p_*,
$$
give
    \begin{equation}
    \label{18}
    \begin{split}
    \left|\int_{\partial\Omega}\mathcal C(x,u_n)(u_n-u)d\sigma \right| &\leq
\int_{\partial\Omega}(c_1 \vert u_n \vert^{\alpha_3}+ c_2)\vert u_n-u \vert d\sigma\\
    &\leq c_1 \|u_n\|_{L^{p_*}(\partial\Omega)}^{\alpha_3}\|u_n-u\|_{L^{\frac{p_*}{p_*-\alpha_3}}(\partial\Omega)} \\
    & \quad + c_2\vert
\partial\Omega\vert^{\frac{p-1}{p}}\|u_n-u\|_{L^p(\partial\Omega)}\to 0\ \quad \mbox{as $n\to\infty$}.
    \end{split}
    \end{equation}
If we gather \eqref{15}, \eqref{16} and \eqref{18}, in  view of
\eqref{operator T} then inequality \eqref{12} becomes
    \begin{equation}\label{12'}
    \limsup_{n\to+\infty} \int_{\Omega} \mathcal A(x,u_n, \nabla u_n)(u_n-u)dx\leq 0.
    \end{equation}
Thanks to assumptions (A1)-(A3) it is allowed to invoke
 \cite[Theorem 2.109]{CLM}. Then \eqref{12'} and the weak convergence
$u_n\rightharpoonup u$ in $W^{1,p}(\Omega)$ ensure the strong
convergence $u_n\to u$ in  $W^{1,p}(\Omega)$. Once the strong
convergence is achieved, it is straightforward to deduce from the
continuity of the involved Nemytskii maps that the nonlinear
operator $T$ is pseudomonotone.

The next step  is to show that  $T$ is coercive. To this end, first
observe that
    \begin{equation}
    \label{coerc}
    \begin{split}
    \langle Tv,v\rangle&=\int_{\Omega} \mathcal A(x,v,\nabla
v)\cdot \nabla v dx+a\int_{\Omega}\vert v\vert^p dx \\
    & \quad +\int_{\Omega}\mathcal B(x,\rho\ast E(v),\nabla(\rho \ast E(v)))
v dx+\int_{\partial\Omega}\mathcal C(x,v)v d\sigma.
    \end{split}
    \end{equation}
We estimate the terms of the inequality above separately. First of
all thanks to assumption (A3)  we have
    \[
    \begin{split}
    \into \mathcal A(x, v, \nabla v) \cdot \nabla v dx&\geq \into (a_4 \vert \nabla v \vert^p- a_5) dx \\
    &= a_4 \Vert \nabla v \Vert_{L^p(\Omega)}^p- a_5 |\Omega|.
    \end{split}
    \]
Moreover, reasoning as in \eqref{bound-f}-\eqref{conv2} we have
    \[
    \begin{split}
    \into \mathcal B(x, \rho \ast E(v), \nabla(\rho \ast E(v))) v & \ge - \into \l[f(x) + b_1 \vert \rho\ast E(v) \vert^{\alpha_1}+ b_2\vert \nabla(\rho \ast E(v)) \vert^{\alpha_2}\r]\vert v\vert dx \\
    & \ge - C_{13} \Vert v \Vert_{L^r(\Omega)}- C_{14} \|v\|_{W^{1,p}(\Omega)}^{\alpha_1}
\|v\|_{L^{\frac{p^*}{p^*-\alpha_1}}(\Omega)} \\
& \quad -C_{15}\|v\|_{W^{1,p}(\Omega)}^{\alpha_2}
\|v\|_{L^{\frac{p}{p-\alpha_2}}(\Omega)}
    \end{split}
    \]
as well as
    \[
    \begin{split}
    \int_{\rand} \mathcal C(x, v) v d\sigma & \ge -\int_{\rand} (c_1 \vert v \vert^{\alpha_3}+ c_2) \vert v \vert d\sigma \\
    & \ge -c_1 \|v\|_{L^{\alpha_3+1}(\partial\Omega)}^{\alpha_3+1}- C_{16} \|v\|_{L^{p}(\partial\Omega)}.
    \end{split}
    \]
From \eqref{coerc} we easily derive
    \[
    \begin{split}
    \langle Tv, v \rangle & \ge a_4 \||\nabla v|\|_{L^p(\Omega)}^p\| + a\| v\|_{L^p(\Omega)}^p
    \\
    & \qquad - C_{17}\bigl(\| v\|_{W^{1,p}(\Omega)}^{\alpha_1+1}+\|v\|_{W^{1,p}(\Omega)}^{\alpha_2+1}+\| v\|_{W^{1,p}(\Omega)}^{\alpha_3+1}+\| v\|_{W^{1,p}(\Omega)}+1\bigr)
    \end{split}
    \]
for every $v \in W^{1,p}(\Omega)$. Then by virtue of hypothesis
\eqref{21} we have
\begin{eqnarray*}
\lim\limits_{\|v\|_{W^{1,p}(\Omega)} \rightarrow+ \infty}
\frac{\langle Tv, v\rangle}{\|v\|_{W^{1,p}(\Omega)}}  = +\infty,
\end{eqnarray*}
thus the coercivity of  $T$  ensues. We have already shown that the
nonlinear operator $T$  is bounded, pseudomonotone and coercive.
Consequently, all the requirements of Theorem \ref{TMPMO} are
fulfilled. Therefore,  there exists $u\in W^{1,p}(\Omega)$ verifying
$T(u)=0$. Taking into account \eqref{operator T} it follows that $u$
is a weak solution to problem \eqref{problem}, which completes the
proof.
\section{Proof of Theorem \ref{T2}}\label{S4}

Let $ u\in W^{1,p}(\Omega)$ be a weak solution to \eqref{problem}
for which we can admit that $u\not\equiv 0$. First, we show that
$u\in L^r(\Omega)$ for every $r \in [1,+\infty)$. According to
\eqref{pm} and to the fact that, in the nonlocal terms, the operator $E$ and the convolution with $\rho$ are linear maps, we can suppose that $u \ge 0 $, otherwise we work with $u^+$ and $u^-$. Moreover, throughout the
proof we will denote by $M_i$, $i \in \N$, constants which depend on
the given data and possibly on the solution itself, and we will
specify the dependance when it will be relevant.

Let $h>0 $ and set $ u_h(x):= \min\{u(x), h\}$ for $ x \in \Omega$.
For every number $\kappa>0$, choose  $\varphi= u u_h^{\kappa p}$ as
test function in \eqref{weak-sol}. We note that
$$
\nabla \varphi= u_h^{\kappa p} \nabla u + \kappa p u u_h^{\kappa
p-1} \nabla u_h.
$$
Inserting such a $\varphi$ in \eqref{weak-sol} gives
\begin{align}\label{1n}
\begin{split}
&\int_{\Omega} (\mathcal{A}(x,u,\nabla u)\cdot \nabla u) u_h^{\kappa
p}dx + \kappa p\int_{\Omega}(\mathcal{A}(x,u,\nabla u) \cdot \nabla
u_h) u_h^{\kappa p-1}udx
+a\int_{\Omega}u^p u_h^{\kappa p} dx\\
&=\int_{\Omega}\mathcal B(x, \rho \ast E(u), \nabla(\rho \ast E(u)))
uu_h^{\kappa p} dx + \int_{\partial \Omega} \mathcal{C}(x,u)
uu_h^{\kappa p} d \sigma.
\end{split}
\end{align}
Applying condition {(H2)} yields
\begin{align}\label{2}
\begin{split}
&\int_{\Omega} (\mathcal{A}(x,u,\nabla u) \cdot \nabla u) u_h^{\kappa p}dx\\
&\geq \int_{\Omega} \left[a_4 |\nabla u|^p-a_5u^{p^*}-a_6\right]u_h^{\kappa p} dx\\
& \geq a_4\int_{\Omega} \left|\nabla u\right|^p u_h^{\kappa p}dx
-(a_5+a_6) \int_{\Omega} u^{p^*}u_h^{\kappa p}dx-a_6|\Omega|
\end{split}
\end{align}
and
\begin{align}\label{2b}
\begin{split}
&\int_{\Omega} (\mathcal{A}(x,u,\nabla u) \cdot \nabla u_h) u_h^{\kappa p-1}udx\\
&=\int_{\{x\in\Omega: \, u(x)\leq h\}} (\mathcal{A}(x,u,\nabla u) \cdot \nabla u) u_h^{\kappa p}dx\\
&\geq \int_{\{x\in\Omega: \, u(x)\leq h\}} \left[a_4 |\nabla u|^p-a_5u^{p^*}-a_6\right]u_h^{\kappa p} dx\\
& \geq a_4 \int_{\{x \in \Omega: \, u(x) \le h\}} \left|\nabla u
\right|^p u_h^{\kappa p}dx -(a_5+a_6) \int_{\Omega}
u^{p^*}u_h^{\kappa p}dx-\kappa pa_6|\Omega|.
\end{split}
\end{align}
Note that in the last passage of both \eqref{2} and \eqref{2b} we use the following fact
	\[
	u_h^{\kappa p} \le u^{p^*} u_h^{\kappa p}+ 1.
	\]
Indeed, if $u> 1$, then $u^{p^*}> 1$, which implies that
	\[
	u_h^{\kappa p} \le u^{p^*} u_h^{\kappa p}< u^{p^*} u_h^{\kappa p}+ 1.
	\]
If $u \le 1$, then we refer to the definition of $u_h:= \min\{u(x), h\}$, and again distinguish among two cases.

If $h>1$, then $u_h(x)= u(x) \le 1$, and it follows that
	\[
	u_h^{\kappa p} \le 1< 1+ u^{p^*}  u_h^{\kappa p},
	\]
because $u^{p^*} u_h^{\kappa p}> 0$. 
If $h\le 1$, then $u_h(x)= h \le 1$, and we have again
	\[
	u_h^{\kappa p} \le 1< 1+ u^{p^*} u_h^{\kappa p}.
	\]
By means of condition {(H3)} we have
\begin{equation}\label{3}
\begin{split}
& \int_{\Omega} \mathcal B(x,\rho \ast E(u),\nabla(\rho \ast E(u))) u u_h^{\kappa p} dx \\
& \le \int_{\Omega} \l( f(x) + b_1 \vert \rho \ast E(u)
\vert^{\alpha_1} + b_2 \vert \nabla (\rho \ast E(u))
\vert^{\alpha_2}\r) u u_h^{\kappa p} dx.
\end{split}
\end{equation}
We estimate the terms on the right-hand side of \eqref{3}
separately. First, through H\"older's inequality we have
\begin{equation}\label{3a}
\int_{\Omega} f(x) u u_h^{\kappa p} dx \le \Vert f \Vert_{r'} \l(
\int_{\Omega} (u u_h^{\kappa p})^r dx \r)^{1/r} \le M_1 (1+ \Vert u
u_h^{\kappa} \Vert_{L^{pr}(\Omega)}^p).
\end{equation}
Moreover, we set $ r_1:= \frac{p^*}{p^*-\alpha_1} $ and $ r_2:=
\frac{p}{p- \alpha_2}$. Making use of H\"older's inequality, with an
argument similar as in \eqref{conv1}-\eqref{conv2}, we find that
    \begin{equation}\label{3b}
    \begin{split}
    \int_{\Omega}\vert \rho \ast E(u)\vert^{\alpha_1} u u_h^{\kappa p}dx
& \le \Vert \rho \ast E(u)\Vert_{L^{p^*}(\mathbb{R}^N)}^{\alpha_1}
\Vert u u_h^{\kappa p} \Vert_{L^{r_1}(\Omega)} \\
    & \le M_2 \Vert \rho \ast E(u)
\Vert_{W^{1,p}(\mathbb{R^N})}^{\alpha_1} \Vert u u_h^{\kappa p}
\Vert_{L^{r_1}(\Omega)}\\
    & \le M_3 \Vert \rho
\Vert_{L^1(\mathbb{R}^N)}^{\alpha_1} \Vert u
\Vert_{W^{1,p}(\Omega)}^{\alpha_1}
\Vert u u_h^{\kappa p} \Vert_{L^{r_1}(\Omega)} \\
    & \le M_4 \l(1+ \Vert u u_h^{\kappa} \Vert_{{L^{pr_1}(\Omega)}}^p
\r)
    \end{split}
    \end{equation}
and
    \begin{equation}\label{3c}
    \begin{split}
    \int_{\Omega} \vert \nabla (\rho \ast E(u))\vert^{\alpha_2} u
u_h^{\kappa p} dx & \le M_5 \Vert |\nabla (\rho \ast
E(u))|\Vert_{L^p(\mathbb{R}^N)}^{\alpha_2}
\Vert u u_h^{\kappa p} \Vert_{L^{r_2}(\Omega)} \\
    & \le M_6 \Vert \rho \Vert_{L^1(\mathbb{R}^N)}^{\alpha_2}
 \Vert|\nabla u|\Vert_{L^{p}(\Omega)}^{\alpha_2}
\Vert u u_h^{\kappa p} \Vert_{L^{r_2}(\Omega)} \\
    & \le M_7 \l(1+ \Vert u u_h^{\kappa} \Vert_{L^{p r_2}(\Omega)}^p \r),
    \end{split}
    \end{equation}
where the constants $M_4$ and $M_7$ depend on the solution $u$,
precisely
\begin{equation}\label{0}
M_4=M_4(\Vert u \Vert_{W^{1,p}(\Omega)}) \quad \text{and} \quad
M_7=M_7(\Vert \nabla u \Vert_{L^{p}(\Omega)}).
\end{equation}
Via hypothesis (H4) we estimate
\begin{align}\label{4}
\begin{split}
\int_{\partial \Omega} \mathcal{C}(x,u) uu_h^{\kappa p} d \sigma
& \leq \int_{\rand} \left(c_1u^{p_*-1}+c_2\right) uu_h^{\kappa p}d\sigma\\
& \leq (c_1+c_2) \int_{\partial\Omega} u^{p_*}u_h^{\kappa p}d\sigma
+c_2|\partial \Omega|.
\end{split}
\end{align}
From \eqref{alpha} and the hypothesis on $r$, we see that
\begin{equation}\label{exponents}
\tilde r:= \max\l\{r, r_1, r_2 \r\} < \frac{p^*}{p}.
\end{equation}
Combining \eqref{1n}-\eqref{3c}, \eqref{4}, \eqref{exponents}
results in
\begin{equation}
\label{5}
\begin{split}
& a_4 \l(\int_{\Omega} \vert \nabla u \vert^p u_h^{\kappa p} dx+
\kappa p
\int_{\{x \in \Omega: \, u(x) \le h\}} \vert \nabla u \vert^p u_h^{\kappa p} dx \r) \\
& \le (\kappa p+1)(a_5+ a_6) \into u^{p^*} u_h^{\kappa p} dx+ (c_1+
c_2)
\int_{\partial \Omega} u^{p_*} u_h^{\kappa p} d\sigma \\
& \qquad+ M_8 \Vert u u_h^{\kappa} \Vert_{L^{p\tilde r}(\Omega)}^p+
M_9(\kappa+1),
\end{split}
\end{equation}
with positive constants $M_8$ and $M_9$ independent on $\kappa$.

Notice that
\begin{equation*}
\begin{split}
& \int_{\Omega}\vert \nabla u \vert^p u_h^{\kappa p} dx+ \kappa p
\int_{\{x \in \Omega: \, u(x) \le h\}} \vert \nabla u \vert^p
u_h^{\kappa p} dx
\\
& =\int_{\{x \in \Omega: \, u(x) > h\}} \vert \nabla u \vert^p
u_h^{\kappa p} dx +(\kappa p+1)
\int_{\{x \in \Omega: \, u(x) \le h\}} \vert \nabla u \vert^p u_h^{\kappa p} dx  \\
& \ge \frac{\kappa p+1}{(\kappa+1)^p} \int_{\{x \in \Om: \, u(x)> h\}} |\nabla u|^p u_h^{\kappa p} dx+ (\kappa p+1) \int_{\{x \in \Om: \, u(x) \le h\}} |\nabla u|^p u_h^{\kappa p} dx \\
& \geq \frac{\kappa p+1}{(\kappa+1)^p}\int_{\Omega} \vert \nabla(u
u_h^{\kappa})\vert^p dx,
\end{split}
\end{equation*}
thanks to Bernoulli's inequality $(\kappa +1)^p \ge \kappa p+1$ and to the fact that $(\kappa+1)^p>1$. 
Therefore, \eqref{5} and \eqref{norms} entail
\begin{equation}
\label{7}
\begin{split}
\frac{\kappa p+1}{(\kappa+1)^p} \Vert u u_h^{\kappa}
\Vert_{W^{1,p}(\Omega)}^p & \le \frac{\kappa p+1}{(\kappa+1)^p}
\Vert u u_h^{\kappa} \Vert_{L^{p}(\Omega)}^p
+ M_{10}(\kappa p+1) \int_{\Omega} u^{p^*} u_h^{\kappa p} dx \\
& \qquad+ M_{11} \int_{\partial\Omega} u^{p_*} u_h^{\kappa p}
d\sigma
+ M_8 \Vert u u_h^{\kappa} \Vert_{L^{p\tilde r}(\Omega)}^p+ M_9(\kappa+1) \\
& \le M_{10} (\kappa p+1) \int_{\Omega}u^{p^*} u_h^{\kappa p} dx
+ M_{11} \int_{\partial\Omega} u^{p_*} u_h^{\kappa p} d\sigma \\
& \qquad + M_{12} \l(\frac{\kappa p+1}{(\kappa+1)^p} + 1 \r) \Vert u
u_h^{\kappa} \Vert_{L^{p\tilde r}(\Omega)}^p+ M_9(\kappa+1).
\end{split}
\end{equation}
We now aim to estimate the critical integrals on the right-hans side of \eqref{7}. To this end, we set $A:= u^{p^*-p}$ and $B:= u^{p_*-p}$, and take $\Lambda, \Gamma>0$. Then H\"older's inequality and the Sobolev embedding give
	\begin{equation}
	\label{critic-eq1}
	\begin{split}
	& \into u^{p^*} u_h^{\kappa p} dx \\
	&= \int_{\{x \in \Om : \, A(x) \le \Lambda\}} A (u u_h^{\kappa})^p dx+ \int_{\{x \in \Om: \, A(x)> \Lambda\}} A(u u_h^{\kappa})^p dx \\
	& \le \Lambda \int_{\{x \in \Om: \, A(x) \le \Lambda} (u u_h^{\kappa})^p dx \\
	& \qquad +\l(\int_{\{x \in \Om: \, A(x)> \Lambda\}} A^{\frac{p^*}{p^*-p}}dx \r)^{\frac{p^*-p}{p^*}} \l(\into (u u_h^{\kappa})^{p^*} dx \r)^{\frac{p}{p^*}} \\
	& \le \Lambda \|u u_h^{\kappa} \|_{L^p(\Om)}^p+ \l(\int_{\{x \in \Om: \, A(x)> \Lambda\}} A^{\frac{p^*}{p^*-p}} dx \r)^{\frac{p^*-p}{p^*}} C_{\Om}^p \|u u_h^{\kappa}\|_{W^{1,p}(\Om)}^p
	\end{split}
	\end{equation}
as well as
	\begin{equation}
	\label{critic-eq2}
	\begin{split}
	& \int_{\rand} u^{p_*} u_h^{\kappa p} d\sigma \\
	&= \int_{\{x \in \rand: \, B(x) \le \Gamma\}} B (uu_h^{\kappa})^p d\sigma+ \int_{\{x \in \rand: \, B(x)> \Gamma\}} B (uu_h^{\kappa})^p d\sigma \\
	& \le \Gamma \int_{\{x \in \rand: \, B(x) \le \Gamma\}} (u u_h^{\kappa})^p d\sigma \\
	& \qquad +\l(\int_{\{x \in \rand: \, B(x)> \Gamma\}} B^{\frac{p_*}{p_*-p}} d\sigma \r)^{\frac{p_*-p}{p_*}} \l(\int_{\rand} (uu_h^{\kappa})^{p_*} d\sigma\r)^{\frac{p}{p_*}} \\
	& \le \Gamma \|u u_h^{\kappa}\|_{L^p(\rand)}^p+ \l(\int_{\{x \in \rand: \, B(x)> \Gamma\}} B^{\frac{p_*}{p_*-p}} d\sigma\r)^{\frac{p_*-p}{p_*}} c_{\rand}^p \|u u_h^{\kappa}\|_{W^{1,p}(\Om)}^p,
	\end{split}
	\end{equation}
with the embedding constants $C_{\Om}$ and $c_{\rand}$. Moreover, if we set
	\begin{equation}
	\label{critic-eq3}
	\begin{split}
	f_1(\Lambda)&:= \l(\int_{\{x \in \Om: \, A(x)> \Lambda\}} A^{\frac{p^*}{p^*-p}} dx \r)^{\frac{p^*-p}{p^*}} \\
	 \text{as well as} \quad f_2(\Gamma)&:= \l(\int_{\{x \in \rand: \, B(x)> \Gamma\}} B^{\frac{p_*}{p_*-p}} d\sigma\r)^{\frac{p_*-p}{p_*}},
	 \end{split}
	\end{equation}
we see that
	\begin{equation}
	\label{critic-eq4}
	f_1(\Lambda) \to 0 \quad \text{as } \Lambda \to 0 \quad \text{as well as} \quad f_2(\Gamma) \to 0 \quad \text{as } \Gamma \to 0.
	\end{equation}
From \eqref{7}, taking into account \eqref{critic-eq1}-\eqref{critic-eq3} and applying H\"older's inequality we have
	\begin{equation}
	\label{7.1}
	\begin{split}
	& \frac{\kappa p+1}{(\kappa+1)^p} \|u u_h^{\kappa}\|_{W^{1,p}(\Om)}^p \\
	& \le M_{13} \l((\kappa p+ 1)\Lambda+ 1+ \frac{\kappa p+1}{(\kappa+1)^p} \r) \|u u_h^{\kappa}\|_{L^{p\tilde r}(\Om)}^p \\
	& \qquad M_{10} (\kappa p+1) f_1(\Lambda) C_{\Om}^p \|u u_h^{\kappa}\|_{W^{1,p}(\Om)}^p+ M_{11} \Gamma \|u u_h^{\kappa}\|_{L^p(\rand)}^p \\
	& \qquad+ M_{11} f_2(\Gamma) c_{\rand}^p \|u u_h^{\kappa}\|_{W^{1,p}(\Om)}^p+ M_9(\kappa+1).
	\end{split}
	\end{equation}
Taking into account \eqref{critic-eq4} we can choose $\Lambda= \Lambda(\kappa, u), \Gamma= \Gamma(\kappa, u)>0$ large enough in order to have
	\[
	M_{10} (\kappa p+1) f_1(\Lambda) C_{\Om}^p= \frac{\kappa p+1}{4(\kappa+1)^p} \quad \text{as well as} \quad M_{11} f_2(\Gamma) c_{\rand}^p= \frac{\kappa p+1}{4(\kappa+1)^p}.
	\]
Then from \eqref{7.1} we have
	\[
	\begin{split}
	& \frac{\kappa p+1}{4(\kappa+1)^p} \|u u_h^{\kappa}\|_{W^{1,p}(\Om)}^p \\
	&  \le M_{13} \l((\kappa p+1) \Lambda(\kappa, u)+ 1+ \frac{\kappa p+1}{(\kappa+1)^p} \r) \|u u_h^{\kappa}\|_{L^{p \tilde r}(\Om)}^p \\
	& \qquad + M_{11} \Gamma(\kappa, u) \|u u_h^{\kappa}\|_{L^p(\rand)}^p+ M_9(\kappa+1),
	\end{split}
	\]
where both $\Lambda(\kappa, u), \Gamma(\kappa, u)$ depend on $\kappa $ and on the solution itself.

From
this point we  proceed as in \cite[Theorem 3.1, Case I.1]{MW} with
 $\Vert u u_h^{\kappa}\Vert_{L^p(\Omega)}$ replaced by $\Vert u
u_h^{\kappa}\Vert_{L^{p\tilde r}(\Omega)}$, which gives us
\begin{align*}\label{20}
\Vert u \Vert_{L^{(\kappa+1)p^*}(\Omega)}\le M_{14}(\kappa, u)
\end{align*}
for any $ \kappa>0$, where $ M_{14}(\kappa, u) $ is a positive
constant which depends on $ \kappa $ and on the solution $u$.
Consequently, the claim that $u\in L^r(\Omega)$ for every $r \in
[1,\infty)$ follows.

Once the $L^r(\Omega)$-bound is reached, the proof of the
$L^r(\rand)$-boundedness is straightforward (see \cite[Case
I.2]{MW}).

We are now in a position to establish the $L^{\infty}$-boundedness
of $u$. Taking advantage of \eqref{exponents}, we fix
 $ q_1 \in (p\tilde r, p^*)$ and $q_2 \in (p, p_*)$. By H\"older's inequality and
the obtained $ L^r$-bounds in $\Omega$ and on $\partial\Omega$, we
can express \eqref{7} in the form
\begin{equation*}\label{29}
\begin{split}
\frac{\kappa p+1}{(\kappa+1)^p}\Vert u u_h^{\kappa}
\Vert_{W^{1,p}(\Omega)}^p &\le M_{15}\l(\frac{\kappa
p+1}{(\kappa+1)^p}
+ \kappa p+ 2\r) \Vert u u_h^{\kappa} \Vert_{L^{q_1}(\Omega)}^p \\
& \qquad + M_{16} \Vert u u_h^{\kappa}
\Vert_{L^{q_2}(\partial\Omega)}^p+ M_{17}(\kappa+1).
\end{split}
\end{equation*}
Then, proceeding as in \cite[Case II.1]{MW}, arranging the constants
and applying H\"older's inequality, the Sobolev embedding and 
Fatou's lemma we achieve
    \[
    \Vert u \Vert_{L^{(\kappa_n+1)p^*}(\Omega)} \le M_{18},
    \]
where $M_{18}$ is independent on $\kappa $ and  $(\kappa_n+ 1)p^*
\to \infty $ as $n \to \infty$.

Therefore, we can invoke Proposition \ref{prop1}, whence $u \in
L^{\infty}(\Omega)$. Finally, by Proposition \ref{prop3}, it follows
that $\gamma u\in L^{\infty}(\partial\Omega)$. The proof is thus
complete.

\begin{remark}\label{R1}
Hypothesis \emph{(H1)} is not needed in the proof of Theorem
\ref{T2}, but it is necessary in order to have a well-defined weak
solution as formulated in \eqref{weak-sol}.
\end{remark}

\begin{remark}\label{R2}
The bounds obtained in Theorems \ref{T2} depend on the data in
assumption \emph{(H)} and on the solution itself. The proof shows
that the following estimate is valid
\begin{align}\label{01}
\|u\|_{L^r(\Omega)} \leq M(\|u\|_{L^{p^*}(\Omega)}), \quad \forall
\, r\geq 1,
\end{align}
with a constant $M(\|u\|_{L^{p^*}(\Omega)})$ depending on
$\|u\|_{p^*}$. The key step for proving estimate \eqref{01} is
\eqref{0}.
\end{remark}

\begin{remark}\label{R3}
Once \eqref{01} is reached, an alternative reasoning to get the
uniform boundedness of $u$ can be carried out as follows. Let
$0<t<\|u\|_{L^{\infty}(\Omega)}$, where a priori one can have
$\|u\|_{L^{\infty}(\Omega)}=+\infty$. Setting
$$
\Omega_t=\{x\in \Omega: |u(x)|>t\},
$$
it is clear that
$$
\|u\|_{L^r(\Omega)}\geq \left(\int_{\Omega_t}|u|^r dx\right
)^{\frac{1}{r}}\geq t|\Omega_t|^{\frac{1}{r}},\quad \forall \, r\geq
1,
$$
so
$$
\liminf_{r\to\infty}\|u\|_{L^r(\Omega)}\geq t.
$$
Since $t\in (0,\|u\|_{\infty})$ is arbitrary, we deduce that
$$
\liminf_{r\to\infty}\|u\|_{L^r(\Omega)} \geq
\|u\|_{L^{\infty}(\Omega)}.
$$
In view of estimate \eqref{01}, the conclusion that $u\in
L^{\infty}(\Omega)$ is achieved.
\end{remark}

\section*{Acknowledgements}
The authors thank the referees for their useful comments that helped
to improve the paper.

\end{document}